\documentclass{article}
\usepackage[intlimits]{amsmath}
\usepackage{amstext, amssymb, amsthm}
\usepackage[mathscr]{euscript}
\usepackage{eufrak}
\usepackage[mathscr]{euscript}
\usepackage{cancel,enumerate}
\usepackage[affil-it]{authblk}
\usepackage{mathtools}

\theoremstyle{definition}

\newtheorem{example}{Example}

\title{Numerical instability in B\&D models}

\date{\today}

\author[1]{Jan Vrbik}
\affil[1]{Department of Mathematics and Statistics\\ Brock University, Canada}

\begin{document}

\maketitle

\begin{abstract}
When computing the expected value till extinction of a Birth and Death process, the usual textbook approach results in an extreme case of numerical ill-conditioning, which prevents us from getting accurate answers beyond the first few low-lying states; in this brief note we present a potential solution. We also present a novel derivation of the related formulas.\\

\noindent\textbf{MCS}: 60J80, 65Q30
\end{abstract}

\section{Probability of extinction}

Consider a Birth \& Death process with rates of $\lambda_{n}$ (for State $n$ to State $n+1$ transition) and $\mu_{n}$ (State $n$ to State $n-1$). When $\lambda_{0}=\mu_{n}=0$, State $0$ is obviously absorbing and the main task is to establish the probability of the process' ultimate extinction (denoted $a_{i}$), given it starts in State $i.$ Based on what happens during the next transition, the sequence of these probabilities meets the following infinite set of linear equations
\begin{equation}
a_{i}=\frac{\lambda_{i}}{\lambda_{i}+\mu_{i}}~a_{i+1}+\frac{\mu_{i}}{\lambda_{i}+\mu_{i}}~a_{i-1}  \label{0}
\end{equation}
where $i=1,\, 2,\, 3,\ldots$ and $a_{0}=1$. The last equation can be rewritten as
\begin{equation}
\frac{\lambda_{i}}{\lambda_{i}+\mu_{i}}(a_{i}-a_{i+1})=\frac{\mu_{i}}{\lambda_{i}+\mu_{i}}(a_{i-1}-a_{i})
\end{equation}
or, equivalently, introducing $d_{i}=a_{i-1}-a_{i}$
\begin{equation}
d_{i+1}=\frac{\mu_{i}}{\lambda_{i}}~d_{i}  \label{di}
\end{equation}
Note that 
\begin{equation}
a_{i}=1-\sum_{k=1}^{i}d_{k}  \label{ai}
\end{equation}

To find the solution (see \cite{math1}) requires to make State $N$ also absorbing, solving the corresponding \emph{finite} set of equations, and
finally taking the $N\to\infty $ limit. Solving (\ref{di}) is trivial and yields
\begin{equation}
d_{i} = d_{1}\prod_{n=1}^{i-1}\frac{\mu_{n}}{\lambda_{n}}
\end{equation}
Since we have made State $N$ is absorbing, $a_{N}=0$ (i.e. State $0$ cannot be reached from State $N$); this can be restated in the following manner 
\begin{equation}
a_{N} 
= 1-\sum_{k=1}^{N}d_{k}
= 1-d_1\sum_{k=1}^{N}\prod_{n=1}^{k-1}\frac{\mu_{n}}{\lambda_{n}}=0
\end{equation}
further implying that
\begin{equation}
d_{1}
=\frac{1}{\sum_{k=1}^{N}\prod_{n=1}^{k-1}\frac{\mu_{n}}{\lambda_{n}}}
\end{equation}
The solution for the $d_{i}$ sequence is thus
\begin{equation}
d_{i}
=
\frac
  {\displaystyle\prod_{n=1}^{i-1}\frac{\mu_{n}}{\lambda_{n}}}
  {\displaystyle\sum_{k=1}^{N}\displaystyle\prod_{n=1}^{k-1}\frac{\mu_{n}}{\lambda_{n}}}
\xrightarrow[N\to\infty]{}
\frac
  {\displaystyle\prod_{n=1}^{i-1}\frac{\mu_{n}}{\lambda_{n}}}
  {\displaystyle\sum_{k=1}^{\infty }\displaystyle\prod_{n=1}^{k-1}\frac{\mu_{n}}{\lambda_{n}}}  \label{fin}
\end{equation}
which, together with (\ref{ai}), yields any of the probabilities of ultimate extinction. Note that (\ref{fin}) and (\ref{ai}) imply that 
$\lim_{i\to \infty }a_{i}=0$ 
when the infinite sum in the last denominator converges, and $a_{i}=1$ for all $i$ when the sum diverges (extinction is then certain for all initial states). Also note that an empty product, i.e. $\prod_{n=1}^{0}\frac{\mu_{n}}{\lambda_{n}}$, equals to $1$.

An alternate approach (see \cite{Karlin}) is to start with 
\begin{align}
a_0   &= 1 \\ 
a_{1} &= 1-\left( \sum_{k=1}^{\infty }\prod_{n=1}^{k-1}\frac{\mu_{n}}{\lambda_{n}}\right) ^{-1}
\end{align}
and then iterate using the following recursive formula 
\begin{equation}
a_{i+1}=\left( 1+\frac{\mu_{i}}{\lambda_{i}}\right) a_{i}-\frac{\mu_{i}}{\lambda_{i}}a_{i-1}
\end{equation}

Both algorithms usually work quite well and yield practically identical results, but our recommendation is to use the former one, namely the direct evaluation of (\ref{fin}) followed by (\ref{ai}).

\section{Expected time till extinction}

When ultimate extinction is certain, the next issue to investigate is: how long does it take to reach State $0$? To find the distribution of this
random variable is too difficult to even attempt; we usually settle for the corresponding expected value (denoted $\omega_{i}$, when starting in State $i$). The solution is found by solving the following analog of (\ref{0}) which similarly relates three consecutive terms of the corresponding sequence of these expected values except now we have to add the expected time till next \emph{transition} to the RHS, getting
\begin{equation}
\omega_{i}
= \frac{\lambda_{i}}{\lambda_{i}+\mu_{i}}~\omega_{i+1}+\frac{\mu_{i}}{\lambda_{i}+\mu_{i}}~\omega_{i-1}+\frac{1}{\lambda_{i}+\mu_{i}}  \label{rec}
\end{equation}
or, equivalently
\begin{equation}
\omega_{i+1}
= \left( 1+\frac{\mu_{i}}{\lambda_{i}}\right) \omega_{i}-\frac{\mu_{i}}{\lambda_{i}}\omega_{i-1}-\frac{1}{\lambda_{i}}  \label{1}
\end{equation}

To solve (\ref{1}), one introduces $\delta_{i}=\omega_{i+1}-\omega_{i}$
which simplifies the previous equation to
\begin{equation}
\delta_{i}
= \frac{\mu_{i}}{\lambda_{i}}\delta_{i-1}-\frac{1}{\lambda_{i}}  \label{2}
\end{equation}
Its formal solution is
\begin{equation}
\delta_{i} = \sum_{n=i+1}^{\infty }\frac{1}{\lambda_{n}}\prod_{j=i+1}^{n}\frac{\lambda_{j}}{\mu_{j}}+c\prod_{j=1}^{i}\frac{\mu_{j}}{\lambda_{j}}  \label{3}
\end{equation}
where $c$ is an arbitrary constant.

\begin{proof}
First we show that the first term of RHS of (\ref{3}) is a particular solution to (\ref{2}):
\begin{align}
\MoveEqLeft\sum_{n=i+1}^{\infty }\frac{1}{\lambda_{n}}\prod\limits_{j=i+1}^{n}\frac{\lambda_{j}}{\mu_{j}}-\frac{\mu_{i}}{\lambda_{i}}\sum_{n=i}^{\infty}\frac{1}{\lambda_{n}}\prod\limits_{j=i}^{n}\frac{\lambda_{j}}{\mu_{j}} \\
&= \sum_{n=i+1}^{\infty }\frac{1}{\lambda_{n}}\prod\limits_{j=i+1}^{n}\frac{\lambda_{j}}{\mu_{j}}-\sum_{n=i}^{\infty }\frac{1}{\lambda_{n}}\prod\limits_{j=i+1}^{n}\frac{\lambda_{j}}{\mu_{j}} \\
&=-\frac{1}{\lambda_{i}}
\end{align}
The second term is clearly a general solution to the homogeneous version of (\ref{2}).
\end{proof}

Note that $\delta_{i}$ can be interpreted as the expected time till reaching State $i$ (for the first time), given the process starts in State $i+1$; this clearly implies that $\delta_{i}$ cannot be a function of any of the $\lambda_{j}$ and $\mu_{j}$ rates when $j$ is less than or equal to $i$, further implying that $c=0$ (note that the first term of (\ref{3}) is already free of the irrelevant rates). The final answer is thus
\begin{equation}
\delta_{i}
= \sum_{n=i+1}^{\infty }\frac{1}{\lambda_{n}}\prod\limits_{j=i+1}^{n}\frac{\lambda_{j}}{\mu_{j}}  \label{fin2}
\end{equation}
with
\begin{equation}
\omega_{i}
= \sum_{k=0}^{i-1}\delta_{k}  \label{omega}
\end{equation}

An alternate approach (used by both \cite{Karlin} and \cite{vrbik}) would be to first evaluate 
\begin{equation}
\omega_{1}
= \delta_{0}
= \sum_{n=1}^{\infty }\frac{1}{\lambda_{n}}\prod\limits_{j=1}^{n}\frac{\lambda_{j}}{\mu_{j}}
\end{equation}
and then use (\ref{1}) recursively to compute as many terms of the $\omega_{i}$ sequence as needed ($\omega_{0}$ is of course equal to $0$). Unlike in the computation of $a_{i}$, this alternate algorithm results in
inaccurate, then incorrect and eventually totally nonsensical values of $\omega_{i}$ as $i$ increases; trying to alleviate this by substantially increasing the accuracy of the computation will only defer the problem to
somehow higher values of $i$.

\begin{figure}[htbp]
\begin{center}
\includegraphics{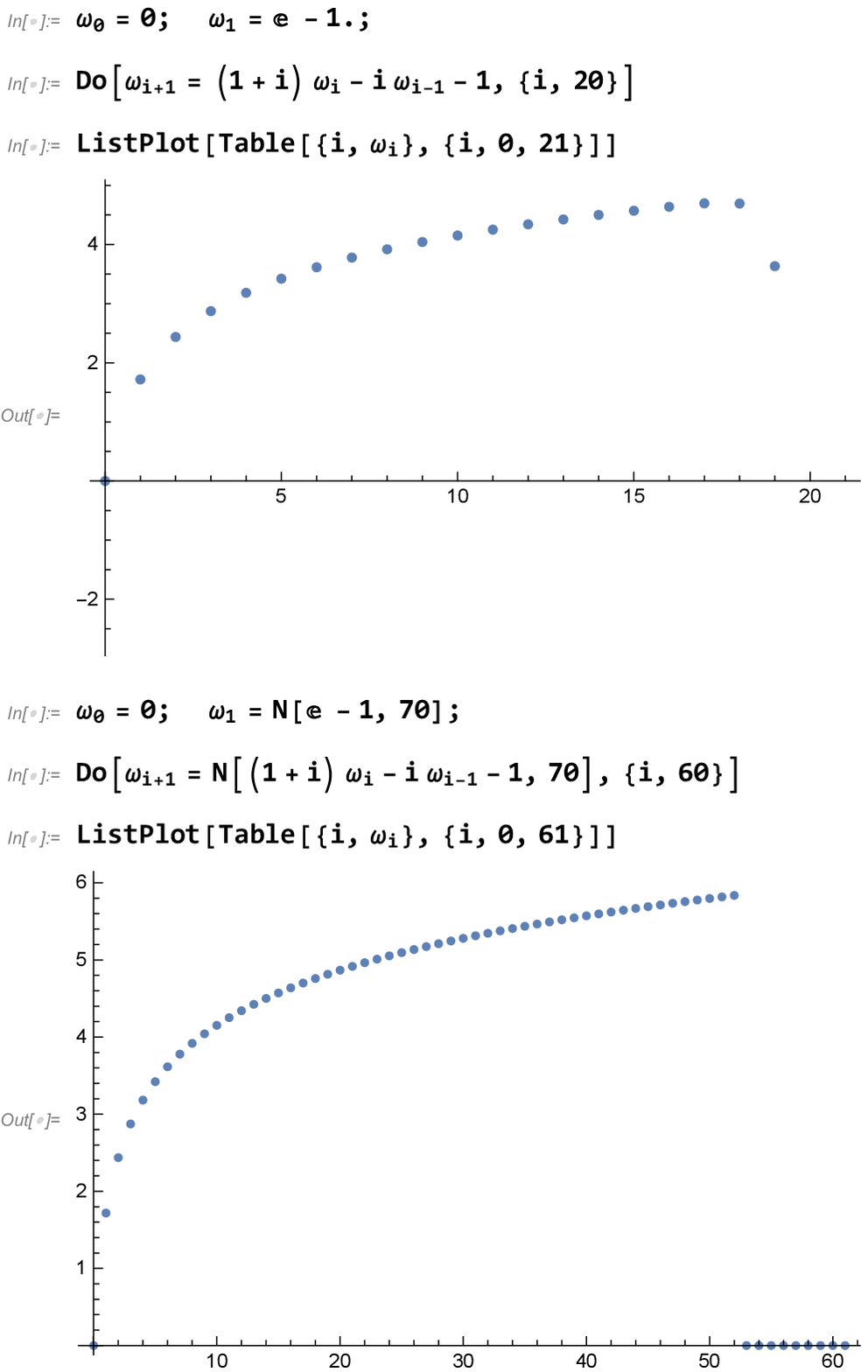}
\caption{Example 1 results.}
\label{math1}
\end{center}
\end{figure}

\begin{example}
Using $\lambda_{n}=1$ and $\mu_{n}=n$ (one of the simplest such models which, furthermore, leads to the following analytic solution: $\omega_{1}=\delta_{0}=e-1$. Using (\ref{1}), repeatedly, to find $\omega_{2}$, $\omega_{3}$, \ldots yields nonsensical results starting at $\omega_{20}$; increasing the accuracy to 70 decimal digits still leads to a similar breakdown at $\omega_{52}$, as the corresponding Mathematica program in Figure \ref{math1} indicates.
\end{example}

Our conclusion (and the main point of this article) is that this algorithm is so badly ill-conditioned that it should never be used.

Instead, we recommend using (\ref{fin2}), followed by (\ref{omega}). But even then, similar numerical ill-conditioning may (rather surprisingly) happen when the right hand side of (\ref{fin2}) is evaluated analytically
(which is possible in some cases). The problem disappears as soon as we switch to numerical evaluation of the same (in Mathematica, this requires using `NSum' instead of `Sum', as the following example demonstrates).

\begin{figure}[htbp]
\begin{center}
\includegraphics{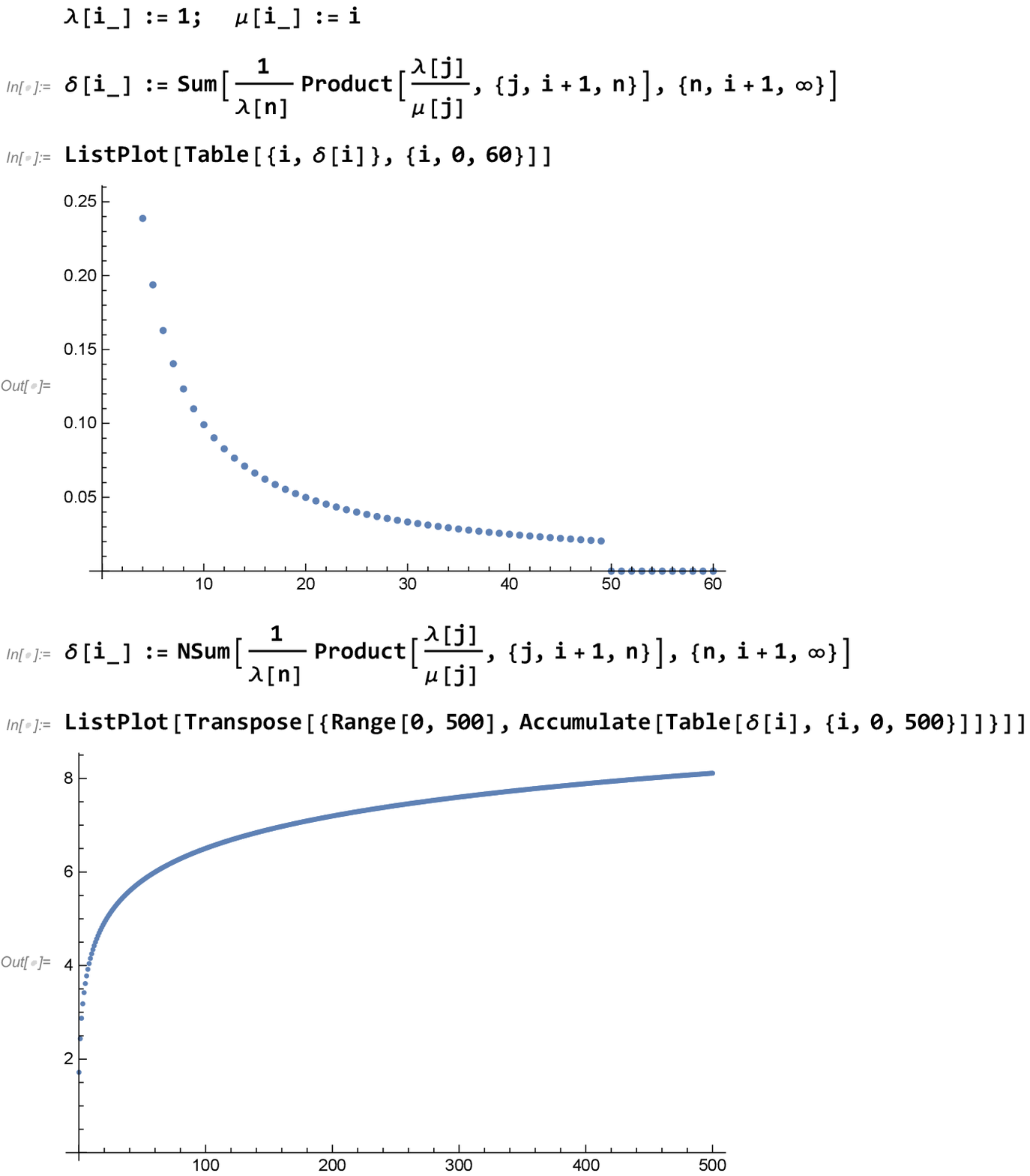}
\caption{Example 2 results.}
\label{math2}
\end{center}
\end{figure}

\begin{example}
We use the same rates as before, but make the code more general. When we allow Mathematica to convert the RHS of (\ref{fin2}) into a formula, its evaluation runs into difficulty at $\omega_{51};$ when we force Mathematica to evaluate the same RHS numerically, correct results are then produced for practically any $\omega_{i}$ (we are demonstrating this up to $\omega_{500}$; in this case we don't display the corresponding $\delta_{i}$ sequence).  See Figure \ref{math2}.
\end{example}

Note that no attempt has bee made to optimize our code -- the computation is fast enough regardless.


\begin{thebibliography}{9}
\bibitem{Karlin} Samuel Karlin: ``A First Course in Stochastic Processes'',
Academic Press, 1969

\bibitem{vrbik} Jan Vrbik and Paul Vrbik: ``Informal Introduction to Stochastic Processes with Maple'', Springer, 2013

\bibitem{math1} Iosif Pinelis: ``Computing probability of ultimate absorption in B\&D processes'', Available: \texttt{https://mathoverflow.net/q/344708} [2021, July]

\bibitem{math2} Iosif Pinelis: ``Expected time till extinction in a B\&D process'', Available: \texttt{https://mathoverflow.net/q/345310} [2021, July]
\end{thebibliography}
\end{document}